# Integral Inequalities for the Analysis of Distributed Parameter Systems: A Complete Characterization via the Least-Squares Principle

Qian Feng[1] and Alexandre Seuret[3] and Sing Kiong Nguang[†]. and Feng Xiao[1,2]

*Abstract*— **A wide variety of integral inequalities (IIs) have been developed and studied for the stability analysis of distributed parameter systems using the Lyapunov functional approach. However, no unified mathematical framework has been proposed that could characterize the similarity and connection between these IIs, as most of them was introduced in a dispersed manner for the analysis of specific types of systems. Additionally, the extent to which the generality of these IIs can be expanded and the optimality of their lower bounds (LBs) remains open questions. In this study, we introduce two general classes of IIs that can generalize nearly all IIs in the literature. The integral kernels of the LBs of our IIs can contain an unlimited number of weighted $\mathscr{L}^2$ functions that are linearly independent in a Lebesgue sense. Moreover, we not only estalish the equivalence relations between the LBs of our IIs, but also demonstrate that these LBs are guaranteed by the least squares principle, implying asymptotic convergence to the upper bound when the kernels functions constitutes a Schauder basis of the underlying Hilbert space. Owing to their general structures, our IIs are applicable in a variety of contexts, such as the stability analysis of coupled PDE-ODE systems or cybernetic systems with delay structures.**

## I. Introduction

A wide variety of control and optimization problems involve the applications of integral inequalities (IIs). Notable examples can be found in the stability analysis and stabilization of linear delay systems [1]–[3] and PDE-related systems [4]–[8], where IIs are crucial for constructing tractable conditions guaranteeing the existence of Lyapunov-Krasovskiĭ functionals (LKFs). To reduce the conservatism induced from the inequality lower bounds (LBs), researchers have developed numerous IIs in various forms over the past decades [1]–[3], [9]–[17]. These developments have directly influenced the structures of LKFs that can be used, which, in turn, affect the conservatism in the stability/synthesis conditions derived from them. Consequently, optimizing the LBs of IIs is essential, as these bounds are vital for effectively applying the LKF approach in the analysis of various distributed parameter systems.

In general, IIs for the analysis of distributed parameter systems can generally be divided into two categories based on their mathematical structures: those without slack variables and those with slack variables. The former follows the research trajectories established in [1]–[3], [10]–[14], [18], where the Jensen [9] and Wirtinger [18] IIs can be regarded as their progenitor. These IIs contain no additional matrices other than the original matrix in the integral quadratic form. In contrast, the latter category of IIs contains free matrix variables, an idea originating from the studies in [19] concerning the stability analysis of systems with time-varying delay. The results established in [20] can be regarded as the first work attempting to provide a unified characterization for these two types of IIs with sufficient generality for the free-matrix type. Nevertheless, the optimality and asymptotic convergence of the LBs were not investigated, and some presentations and notations in [20] are not that accessible for non-professional mathematicians. Consequently, there is a noticeable lack of a unified framework in the literature that can expound the relationship between most of the existing IIs, their limitations, and the optimality of their LBs, especially in more general terms.

In this note, we propose two general classes of IIs following the ideas in [1]–[3], [10]–[14], [20] and [15]–[17], where IIs were gradually generalized to construct Lyapunov-Krasovskiĭ functionals with increasingly general structures. Our IIs can generalize almost all existing results in the literature concerning the stability analysis of linear distributed parameter systems such as time-delay [1]–[3], [10], [11] and PDE related systems [4]–[8]. The first type of our IIs features an upper bound quadratic type integral with a weighted function, where its LB has no extra matrix variables other than the original matrix term in the upper bound. Unlike the existing IIs in [1]–[3], [10]–[14], where the integral kernels of LBs or weight functions were restricted to a specific class of functions, the kernels of our IIs can be any function belonging to a weighted $\mathscr{L}^2$ space with a general weight function. This also means that many existing IIs in [1]–[3], [10]–[14] are the particular cases of the proposed IIs. Moreover, we formally establish in this note that the optimality of the LBs of our IIs can be ensured by the least squares principle in Hilbert space [21], where the kernel functions in the LB do not need to be mutually orthogonal. More

This work was partially supported by the ANR (France) Project ANR-15-CE23-0007 and ANR-18-CE40-0022-01, MCIN/AEI and FEDER (Spain), Grant Numbers: PID2019-109071RB-I00, PID2019-105890RJ-I00, the National Natural Science Foundation of China under Grant Nos. 62303180 and 62273145, Beijing Natural Science Foundation under Grant No. 4222053, and Fundamental Research Funds for Central Universities under Grant 2023MS032, China.

1. School of Control and Computer Engineering, North China Electric Power University Beijing, China. Email: qianfeng@ncepu.edu.cn, qfen204@aucklanduni.ac.nz

2. State Key Laboratory of Alternate Electrical Power System with Renewable Energy Sources, North China Electric Power University, Beijing, China. Emails: fengxiao@ncepu.edu.cn

3. Dpto de Ingeniería de Sistemas y Automática, Universidad de Sevilla, Camino de los Descubrimientos, Sevilla, Spain. Emails: aseuret@us.es

†. Department of Electrical, Computer, and Software Engineering, The University of Auckland, New Zealand. Professor Sing Kiong Nguang passed away in 2020.

importantly, we demonstrate that asymptotic convergence to the upper bound can be established if the kernels of the LBs are a Schauder basis of the underlying Hilbert space.

Next, the second class of IIs, called Free Matrix Type (FMT), is introduced in Section III, which can generalize the existing IIs in [15]–[17]. This type of IIs can be particularly useful when addressing time-vary delays or sampled-data systems as noted in Remark 7 of [16] and Section IV. B of [22], respectively. We then establish a significant finding: under appropriate prerequisites, the LBs of the first and second classes of IIs are equivalent. This discovery implies that the LBs of many existing IIs, previously considered distinct, are actually equivalent. Note that the contents related to FMTII is omitted in this note due to limited space, which will be presented in the journal version.

The rest of the paper is organized as follows. Section II contains the derivation of the first class of IIs, and illustrates its generality using the IIs from existing literature. Additionally, this section provides a rigorous proof of the optimality and asymptotic convergence of the LBs utilizing the least squares principle [21]. We then introduce our FMTII in Section III followed by the paper's conclusion. Please note that the contents starting from Section III are omitted, and will be presented in the journal version of this paper.

*Notation*

Let $\mathbb{S}^n := \{X \in \mathbb{R}^{n\times n} : X = X^\top\}$ and $\mathbb{R}_{\geq a} := \{x \in \mathbb{R} : x \geq a\}$ and $\mathbb{R}^{n\times m}_{[r]} = \{X \in \mathbb{R}^{n\times m} : \mathrm{rank}(X) = r\}$. $\mathcal{M}(\mathcal{X};\mathbb{R}^d)$ denotes the set containing all measurable functions defined from Lebesgue measurable set $\mathcal{X}$ to $\mathbb{R}^d$ endowed with the Boreal algebra. $\mathsf{Sy}(X) := X + X^\top$ stands for the sum of a matrix with its transpose. $[x_i]_{i=1}^n := \begin{bmatrix} x_1^\top \cdots x_i^\top \cdots x_n^\top \end{bmatrix}^\top$ denotes a column vector containing a sequence of mathematical objects such as scalars, vectors, matrices, etc, whereas the double bracket notation $[\![x_i]\!]_{i=1}^n = \begin{bmatrix} x_1 \cdots x_i \cdots x_n \end{bmatrix}$ indicates the row vector counterpart. The symbol $*$ is used to indicate $[*]YX = X^\top Y X$ or $X^\top Y[*] = X^\top Y X$ or $\left[\begin{smallmatrix} A & B \\ * & C \end{smallmatrix}\right] = \left[\begin{smallmatrix} A & B \\ B^\top & C \end{smallmatrix}\right]$. Moreover, $\mathsf{O}_{n\times n}$ denotes an $n\times n$ zero matrix abbreviated by $\mathsf{O}_n$, while $\mathbf{0}_n$ represents an $n\times 1$ column vector. We frequently use $X \oplus Y = \left[\begin{smallmatrix} X & \mathsf{O} \\ * & Y \end{smallmatrix}\right]$ to denote the diagonal sum of two matrices, with its $n$-ary form $\bigoplus_{i=1}^{\nu} X_i = X_1 \oplus X_2 \oplus \cdots \oplus X_\nu$. $\otimes$ stands for the Kronecker product. Finally, we assume the order of matrix operations to be *matrix (scalars) multiplications* $> \otimes > \oplus > +$.

## II. The first integral inequality

In this section, we propose the first class of IIs, and demonstrate its deep connections to the least squares approximations in Hilbert space [21]. We initiate our presentation with several lemmas that are extensively employed throughout the paper.

**Lemma 1.** $\forall X \in \mathbb{R}^{n\times m},\ \forall Y \in \mathbb{R}^{m\times p},\ \ \forall Z \in \mathbb{R}^{q\times r},$

$$(X \otimes I_q)(Y \otimes Z) = (XY) \otimes (I_q Z) = (XY) \otimes Z \\ = (XY) \otimes (Z I_r) = (X \otimes Z)(Y \otimes I_r). \quad (1)$$

The following lemma is partially taken from [22, Lemma 4.1], which is utilized for the derivations in the following sections.

**Lemma 2.** *Let $C \in \mathbb{S}^m_{\succ 0}$ and $B \in \mathbb{R}^{m\times n}$ be given. Then*

$$\forall M \in \mathbb{R}^{m\times n},\ \ B^\top C^{-1} B \succeq M^\top B + B^\top M - M^\top C M \quad (2)$$

*where* (2) *becomes an equality with $M = C^{-1}B$.*

Let us first define the weighted Lebesgue function space

$$\mathscr{L}^2_{\varpi}(\mathcal{K};\mathbb{R}^d) := \left\{\boldsymbol{\phi}(\cdot) \in \mathcal{M}(\mathcal{K};\mathbb{R}^d) : \|\boldsymbol{\phi}(\cdot)\|_{2,\varpi} < \infty\right\} \quad (3)$$

with semi-norm $\|\boldsymbol{\phi}(\cdot)\|^2_{2,\varpi} := \int_{\mathcal{K}} \varpi(\tau)\boldsymbol{\phi}^\top(\tau)\boldsymbol{\phi}(\tau)\mathsf{d}\tau$, where $\varpi(\cdot) \in \mathcal{M}(\mathcal{K};\mathbb{R}_{\geq 0})$ and $\varpi(\cdot)$ has only countably infinite or finite numbers of zero values. Furthermore, $\mathcal{K}$ satisfies $\mathcal{K} \subseteq \mathbb{R} \cup \{\pm\infty\}$ and $\int_{\mathcal{K}} \mathsf{d}\tau \neq 0$.

Now we are ready to set out the first class of IIs.

**Theorem 1.** *Let $\varpi(\cdot)$, $\mathcal{K}$ and $d \in \mathbb{N}$ in* (3) *be given, and assume $\boldsymbol{f}(\cdot) = [f_i(\cdot)]_{i=1}^d \in \mathscr{L}^2_{\varpi}(\mathcal{K};\mathbb{R}^d)$ that satisfies the Gramian criteria [23, Theorem 7.2.10]*

$$\int_{\mathcal{K}} \varpi(\tau)\boldsymbol{f}(\tau)\boldsymbol{f}^\top(\tau)\mathsf{d}\tau \succ 0 \quad (4)$$

*where $\succ 0$ indicates positive definiteness in $\mathbb{S}^d$. Then*

$$\int_{\mathcal{K}} \varpi(\tau)\boldsymbol{x}^\top(\tau)U\boldsymbol{x}(\tau)\mathsf{d}\tau \geq [*]\,(\mathsf{F}\otimes U)\int_{\mathcal{K}} \varpi(\tau)F(\tau)\boldsymbol{x}(\tau)\mathsf{d}\tau \quad (5)$$

*holds for all $\boldsymbol{x}(\cdot) \in \mathscr{L}^2_{\varpi}(\mathcal{K};\mathbb{R}^n)$ and $U \in \mathbb{S}^n_{\succ 0}$, where $F(\tau) = \boldsymbol{f}(\tau) \otimes I_n$ and $\mathsf{F}^{-1} = \int_{\mathcal{K}} \varpi(\tau)\boldsymbol{f}(\tau)\boldsymbol{f}^\top(\tau)\mathsf{d}\tau$. Moreover, the optimality of* (5) *is guaranteed by* (2).

*Proof.* First of all, let $\boldsymbol{\varepsilon}(\tau) = \boldsymbol{x}(\tau) - F^\top(\tau)\boldsymbol{\omega}$ with $\boldsymbol{\omega} \in \mathbb{R}^{dn}$ and $F(\tau) = \boldsymbol{f}(\tau) \otimes I_n$, through which we can obtain

$$\int_{\mathcal{K}} \varpi(\tau)\boldsymbol{\varepsilon}^\top(\tau)U\boldsymbol{\varepsilon}(\tau)\mathsf{d}\tau = \int_{\mathcal{K}} \varpi(\tau)\boldsymbol{x}^\top(\tau)U\boldsymbol{x}(\tau)\mathsf{d}\tau \\ - 2\int_{\mathcal{K}} \varpi(\tau)\boldsymbol{x}^\top(\tau)UF^\top(\tau)\mathsf{d}\tau\boldsymbol{\omega} \\ + \boldsymbol{\omega}^\top\int_{\mathcal{K}} \varpi(\tau)F(\tau)UF^\top(\tau)\mathsf{d}\tau\boldsymbol{\omega}. \quad (6)$$

Apply (1) to $UF^\top(\tau)$ with $F(\tau) = \boldsymbol{f}(\tau) \otimes I_n$. Then

$$UF^\top(\tau) = U(\boldsymbol{f}^\top(\tau) \otimes I_n) = \boldsymbol{f}^\top(\tau) \otimes U = F^\top(\tau)(I_d \otimes U) \quad (7)$$

by which $\int_{\mathcal{K}} \varpi(\tau)\boldsymbol{x}^\top(\tau)UF^\top(\tau)\mathsf{d}\tau\boldsymbol{\omega}$ can be rewritten as

$$\int_{\mathcal{K}} \varpi(\tau)\boldsymbol{x}^\top(\tau)UF^\top(\tau)\mathsf{d}\tau\boldsymbol{\omega} \\ = \int_{\mathcal{K}} \varpi(\tau)\boldsymbol{x}^\top(\tau)F^\top(\tau)\mathsf{d}\tau\,(I_d \otimes U)\,\boldsymbol{\omega} = \boldsymbol{\vartheta}^\top(I_d \otimes U)\boldsymbol{\omega} \quad (8)$$

with $\boldsymbol{\vartheta} = \int_{\mathcal{K}} \varpi(\tau)F(\tau)\boldsymbol{x}(\tau)\mathsf{d}\tau$. On the other hand, we have

$$\int_{\mathcal{K}} \varpi(\tau)F(\tau)UF^\top(\tau)\mathsf{d}\tau = \int_{\mathcal{K}} \varpi(\tau)F(\tau)F^\top(\tau)\mathsf{d}\tau(I_d \otimes U) = \\ \left[\int_{\mathcal{K}} \varpi(\tau)\boldsymbol{f}(\tau)\boldsymbol{f}^\top(\tau)\mathsf{d}\tau \otimes I_n\right](I_d \otimes U) = \mathsf{F}^{-1} \otimes U \quad (9)$$

by (7) and (1), where $\mathsf{F}^{-1} \succ 0$ is defined in the statement of this theorem. Now (6) can be rewritten as

$$\int_{\mathcal{K}} \varpi(\tau)\boldsymbol{\varepsilon}^\top(\tau)U\boldsymbol{\varepsilon}(\tau)\mathsf{d}\tau = \int_{\mathcal{K}} \varpi(\tau)\boldsymbol{x}^\top(\tau)U\boldsymbol{x}(\tau)\mathsf{d}\tau - \mathsf{Sy}\left(\boldsymbol{\vartheta}^\top\left(I_d\otimes U\right)\boldsymbol{\omega}\right) + \boldsymbol{\omega}^\top\left(\mathsf{F}^{-1}\otimes U\right)\boldsymbol{\omega} \geq 0 \quad (10)$$

by (8) and (9) with $U \succ 0$ in (10). Moreover, eq.(10) can be reformulated as

$$\int_{\mathcal{K}} \varpi(\tau)\boldsymbol{x}^\top(\tau)U\boldsymbol{x}(\tau)\mathsf{d}\tau \geq \mathsf{Sy}\left(\boldsymbol{\vartheta}^\top\left(I_d\otimes U\right)\boldsymbol{\omega}\right) - \boldsymbol{\omega}^\top\left(\mathsf{F}^{-1}\otimes U\right)\boldsymbol{\omega} \quad (11)$$

for all $\boldsymbol{x}(\cdot) \in \mathscr{L}^2_{\varpi}(\mathcal{K};\mathbb{R}^n)$ and $U \in \mathbb{S}^n_{\succ 0}$.

Our task now is to construct the LB in (5) via Lemma 2. Applying (2) to (11) with $C = \mathsf{F}^{-1}\otimes U$, $B = (I_d\otimes U)\boldsymbol{\vartheta}$ and $M = \boldsymbol{\omega}$, one can obtain

$$\mathsf{Sy}\left(\boldsymbol{\vartheta}^\top\left(I_d\otimes U\right)\boldsymbol{\omega}\right) - \boldsymbol{\omega}^\top\left(\mathsf{F}^{-1}\otimes U\right)\boldsymbol{\omega} \leq \boldsymbol{\vartheta}^\top(I_d\otimes U)\big(\mathsf{F}\otimes U^{-1}\big)(I_d\otimes U)\,\boldsymbol{\vartheta} = \boldsymbol{\vartheta}^\top(\mathsf{F}\otimes U)\boldsymbol{\vartheta} \quad (12)$$

for all $\boldsymbol{x}(\cdot) \in \mathscr{L}^2_{\varpi}(\mathcal{K};\mathbb{R}^n)$ and $\boldsymbol{\omega} \in \mathbb{R}^{dn}$, where $\leq$ in (12) becomes an equality $=$ with $M = \boldsymbol{\omega} = C^{-1}B = (\mathsf{F}\otimes I_n)\boldsymbol{\vartheta}$. Consequently, (11) gives (5) with $\boldsymbol{\omega} = (\mathsf{F}\otimes I_n)\boldsymbol{\vartheta}$. ∎

Notice that the following II

$$\forall \mathbf{q} \in \mathbb{R}^\nu, \int_{\mathcal{K}} \varpi(\tau)\boldsymbol{x}^\top(\mathbf{q},\tau)U\boldsymbol{x}(\mathbf{q},\tau)\mathsf{d}\tau \geq [*]\left(\mathsf{F}\otimes U\right)\int_{\mathcal{K}} \varpi(\tau)F(\tau)\boldsymbol{x}(\mathbf{q},\tau)\mathsf{d}\tau \quad (13)$$

can be immediately obtained from (5) since $\mathbf{q}$ is independent of the variable of the differential $\mathsf{d}\tau$. The form in (13) is more suitable for the analysis of PDE-related systems [4]–[8] with a spatial variable $\mathbf{q} \in \mathbb{R}^\nu$, $\nu \in \mathbb{N}$. In the subsequent writing, all $\boldsymbol{x}(\mathbf{q},\tau)$ is abbreviated as $\boldsymbol{x}(\tau)$ whenever such an extension is applicable.

*Remark* 1. From [23, Theorem 7.2.10], matrix inequality (4) indicates that all the functions $\{f_i(\cdot)\}_{i=1}^d$ in $\boldsymbol{f}(\cdot)$ are linearly independent in a Lebesgue sense over $\mathcal{K}$. Since (5) holds true for all $\boldsymbol{f}(\cdot) \in \mathscr{L}^2_{\varpi}\big(\mathcal{K};\mathbb{R}^d\big)$ with (4) and appropriate $\varpi(\cdot)$, its generality is evident. For instance, $\boldsymbol{f}(\cdot)$ can contain a list of orthogonal functions as in [3], elementary functions as in [2] or other types of functions as long as they are linearly independent, which can be determined by verifying the inequality in (4) numerically.

### *A. Demonstrating the generality of* (5) *via examples*

In this subsection, we use specific examples from the literature to demonstrate the generality of (5). To do so, let us first consider the standard expression of Jacobi polynomials [24, Section 22.3.2]

$$j_d^{\alpha,\beta}(\tau)_{-1}^1 := \frac{\gamma(d+1+\alpha)}{d!\gamma(d+1+\alpha+\beta)} \times \sum_{k=0}^d \binom{d}{k}\frac{\gamma(d+k+1+\alpha+\beta)}{\gamma(k+1+\alpha)}\left(\frac{\tau-1}{2}\right)^k \quad (14a)$$

over $\tau \in [-1,1]$, where $d \in \mathbb{N}_0$, $\alpha > -1$, $\beta > -1$ and $\gamma(\cdot)$ stands for the standard gamma function. The orthogonality of the standard Jacobi polynomials can be expressed as [24, Section 22.2.1]

$$\int_{-1}^1 (1-\tau)^\alpha(\tau+1)^\beta \boldsymbol{j}_d^{\alpha,\beta}(\tau)_{-1}^1\left[\boldsymbol{j}_d^{\alpha,\beta}(\tau)_{-1}^1\right]^\top \mathsf{d}\tau = \bigoplus_{k=0}^d \frac{2^{\alpha+\beta+1}\gamma(k+\alpha+1)\gamma(k+\beta+1)}{k!(2k+\alpha+\beta+1)\gamma(k+\alpha+\beta+1)} \quad (14b)$$

with $\boldsymbol{j}_d^{\alpha,\beta}(\tau)_{-1}^1 = \left[j_i^{\alpha,\beta}(\tau)_{-1}^1\right]_{i=0}^d$, where $j_i^{\alpha,\beta}(\tau)_{-1}^1$ are orthogonal with respect to $(1-\tau)^\alpha(\tau+1)^\beta$ over $[-1,1]$. On the other hand, it is preferable to derive a general expression for the Jacobi polynomials defined over $[a,b]$ with $b > a$. Specifically, consider the affine substitution $\frac{2\tau-a-b}{b-a} \to \tau$ where $\frac{2\tau-a-b}{b-a}$ satisfies $-1 \leq \frac{2\tau-a-b}{b-a} \leq 1$ for all $\tau \in [a,b]$ with $b > a$. Then $j_d^{\alpha,\beta}\left(\frac{2\tau-a-b}{b-a}\right)_{-1}^1$ can be written as

$$j_d^{\alpha,\beta}(\tau)_a^b := j_d^{\alpha,\beta}\left(\frac{2\tau-a-b}{b-a}\right)_{-1}^1 = \frac{\gamma(d+1+\alpha)}{d!\gamma(d+1+\alpha+\beta)}\sum_{k=0}^d\binom{d}{k}\frac{\gamma(d+k+1+\alpha+\beta)}{\gamma(k+1+\alpha)}\left(\frac{\tau-b}{b-a}\right)^k \quad (15)$$

with $\tau \in [a,b]$. Now using $\frac{2\tau-a-b}{b-a} \to \tau$ to (14b) yields

$$\begin{aligned}&\int_a^b \left(\frac{-2\tau+2b}{b-a}\right)^\alpha\left(\frac{2\tau-2a}{b-a}\right)^\beta \times \boldsymbol{j}_d^{\alpha,\beta}\left(\frac{2\tau-a-b}{b-a}\right)_{-1}^1 [*]\,\mathsf{d}\left(\frac{2\tau-a-b}{b-a}\right) \\ &= \frac{2^{\alpha+\beta+1}}{(b-a)^{\alpha+\beta+1}}\int_a^b (b-\tau)^\alpha(\tau-a)^\beta \boldsymbol{j}_d^{\alpha,\beta}(\tau)_a^b\left[\boldsymbol{j}_d^{\alpha,\beta}(\tau)_a^b\right]^\top\mathsf{d}\tau \\ &= \bigoplus_{k=0}^d \frac{2^{\alpha+\beta+1}\gamma(k+\alpha+1)\gamma(k+\beta+1)}{k!(2k+\alpha+\beta+1)\gamma(k+\alpha+\beta+1)}\end{aligned} \quad (16a)$$

where $\boldsymbol{j}_d^{\alpha,\beta}(\tau)_a^b := \left[j_k^{\alpha,\beta}(\tau)_a^b\right]_{i=0}^d$ for all $\tau \in [a,b]$. Moreover, the equality in (16a) can be rewritten as

$$\mathsf{F}^{-1} = \int_a^b (b-\tau)^\alpha(\tau-a)^\beta \boldsymbol{j}_d^{\alpha,\beta}(\tau)_a^b\left[\boldsymbol{j}_d^{\alpha,\beta}(\tau)_a^b\right]^\top\mathsf{d}\tau = \bigoplus_{k=0}^d \frac{(b-a)^{\alpha+\beta+1}\gamma(k+\alpha+1)\gamma(k+\beta+1)}{k!(2k+\alpha+\beta+1)\gamma(k+\alpha+\beta+1)}, \quad (16b)$$

expressing the orthogonality of (14a) with regards to the weight function $\varpi(\tau) = (b-\tau)^\alpha(\tau-a)^\beta$. Note that with $\alpha = \beta = 0$, the Jacobi polynomials in (15) become the Legendre polynomials [24, Section 22.2.10]

$$j_d^{0,0}(\tau)_a^b = \ell_d(\tau)_a^b := \sum_{k=0}^d\binom{d}{k}\binom{d+k}{k}\left(\frac{\tau-b}{b-a}\right)^k \quad (17)$$

satisfying $\int_a^b \boldsymbol{\ell}_d(\tau)\boldsymbol{\ell}_d^\top(\tau)\mathsf{d}\tau = \bigoplus_{k=0}^d \frac{b-a}{(2k+1)}$ with $d \in \mathbb{N}_0$ and $\tau \in [a,b]$.

In addition to the Legendre polynomials in (17), there are two special cases of (15) that have been frequently investigated in the literature without being recognized as special cases of Jacobi polynomials. These two instances are

$$
\begin{aligned}
&j_d^{\alpha,0}(\tau)_a^b, \ \ \alpha \in \mathbb{N}, \ \ \beta = 0, \ \ \varpi(\tau) = (\tau - a)^{\beta}, \\
&\qquad \mathsf{F}^{-1} = \bigoplus_{k=0}^{d} \frac{(b-a)^{\alpha+1}}{2k+1+\alpha} \qquad (18a)\\
&j_d^{0,\beta}(\tau)_a^b, \ \ \alpha = 0, \ \ \beta \in \mathbb{N}, \ \ \varpi(\tau) = (\tau - a)^{\beta}, \\
&\qquad \mathsf{F}^{-1} = \bigoplus_{k=0}^{d} \frac{(b-a)^{\beta+1}}{2k+1+\beta}. \qquad (18b)
\end{aligned}
$$

In fact, the polynomials in (18a)–(18b) are of particular interest in our study given the Cauchy formula [25, page 193] for repeated integrations:

$$
\begin{aligned}
&\alpha! \int_a^b \int_a^{\tau_1} \cdots \int_a^{\tau_n} f(\tau_{n+1}) \prod_{i=1}^{n+1} \mathsf{d}\tau_i = \int_a^b (b-\tau)^{\alpha} \boldsymbol{x}(t+\tau) \mathsf{d}\tau, \\
&\beta! \int_a^b \int_{\tau_1}^{b} \cdots \int_{\tau_n}^{b} f(\tau_{n+1}) \prod_{i=1}^{n+1} \mathsf{d}\tau_i = \int_a^b (\tau - a)^{\beta} \boldsymbol{x}(t+\tau) \mathsf{d}\tau.
\end{aligned} \qquad (19)
$$

*Remark* 2. The Jacobi polynomials in (18b) are identical to those of orthogonal hyper-geometric polynomials defined in [11, eq.(13)-(14)]. By (19), we can conclude that (18a)–(18b) are equivalent to the polynomials defined in [13, eq.(3)–(4)], [26, eq.(5)-(7)] and [27, eq.(10)].

Having presented the definitions of Jacobi and Legendre polynomials, we furnish a list of existing IIs in Table I as the special cases of (5) with appropriate $\mathcal{K}$, $\varpi(\tau)$, $\boldsymbol{f}(\tau)$ and $\boldsymbol{x}(\cdot)$. It is worth noting that some of the conclusions in Table I are established using the Cauchy formula for repeated integration [25, page 193] in (19).

| (5) | $\mathcal{K}$ | $\varpi(\tau)$ | $\boldsymbol{f}(\tau)$ | $\boldsymbol{x}(\cdot)$ |
|---|---|---|---|---|
| (5) in [1] | $[a,b]$ | 1 | $\boldsymbol{\ell}_d(\tau)$ | $\boldsymbol{x}(\cdot)$ |
| (6) in [1] | $[a,b]$ | 1 | $\boldsymbol{\ell}_d(\tau)$ | $\dot{\boldsymbol{x}}(\cdot)$ |
| (5) in [3] | $\mathcal{O}$ | $w(\tau)$ | $\boldsymbol{m}(\tau)$ | $\boldsymbol{x}(\cdot)$ |
| (16) in [2] | $\mathcal{K}$ | 1 | $\boldsymbol{m}(\tau)$ | $\boldsymbol{x}(\cdot)$ |
| (5) in [28] | $[a,b]$ | 1 | $\begin{bmatrix}1\\ \boldsymbol{p}(\tau)\end{bmatrix}$ | $\boldsymbol{x}(\cdot)$ |
| (27) in [11] | $[a,b]$ | $(\tau-a)^p$ | $\boldsymbol{j}_d^{0,p}(\tau)_a^b$ | $\boldsymbol{x}(\cdot)$ |
| (34) in [11] | $[a,b]$ | $(\tau-a)^p$ | $\boldsymbol{j}_d^{0,p}(\tau)_a^b$ | $\dot{\boldsymbol{x}}(\cdot)$ |
| (1) in [13] | $[a,b]$ | $(\tau-a)^{m-1}$ | $\boldsymbol{j}_d^{0,m-1}(\tau)_a^b$ | $\boldsymbol{x}(\cdot)$ |
| (2) in [13] | $[a,b]$ | $(b-\tau)^{m-1}$ | $\boldsymbol{j}_d^{m-1,0}(\tau)_a^b$ | $\boldsymbol{x}(\cdot)$ |
| (8) in [26] | $[a,b]$ | $(b-\tau)^{m}$ | $\boldsymbol{j}_d^{0,m}(\tau)_a^b$ | $\boldsymbol{x}(\cdot)$ |
| (9) in [26] | $[a,b]$ | $(\tau-a)^{m}$ | $\boldsymbol{j}_d^{m,0}(\tau)_a^b$ | $\boldsymbol{x}(\cdot)$ |
| (2) in [12] | $[a,b]$ | $(\tau-a)^{k}$ | $\boldsymbol{j}_d^{k,0}(\tau)_a^b$ | $\boldsymbol{x}(\cdot)$ |
| (9) in [14] | $[0,+\infty]$ | $K(\tau)$ | $\begin{bmatrix}1\\ g(\tau)\end{bmatrix}$ | $\boldsymbol{x}(\cdot)$ |
| (22) in [29] | $[-h,0]$ | $(-\tau)^b(\tau+h)^a$ | $\boldsymbol{j}_d^{b,a}(\tau)_{-h}^0$ | $\dot{\boldsymbol{x}}(\cdot)$ |

TABLE I: List of IIs as special cases of (5)

Now, let us elaborate on the inclusions of the IIs in Table I. The function $\boldsymbol{p}(\cdot)$ of [28, eq.(5)] in Table I satisfies

$$
\boldsymbol{p}(\cdot) \in \Big\{ \boldsymbol{f}(\cdot) = [f_i(\cdot)]_{i=1}^{d-1} \in \mathscr{L}^2\big([a,b];\mathbb{R}^{d-1}\big) : \int_a^b \boldsymbol{f}(\tau)\boldsymbol{f}^{\top}(\tau)\mathsf{d}\tau = \bigoplus_{i=1}^{d-1} \int_a^b f_i^2(\tau)\mathsf{d}\tau \ \ \& \ \ \int_a^b \boldsymbol{f}(\tau)\mathsf{d}\tau = \mathbf{0}_d \Big\}
$$

corresponding to the auxiliary functions generated by [28, Lemma 1]. Moreover, $\boldsymbol{p}_k(\cdot) = \boldsymbol{j}_d^{k,0}(\tau)_a^b$ in [12, eq.(2)] according to the definition of Jacobi polynomials, as reflected in Table I. The inequality in [2, eq.(16)] corresponds to (5) with $\varpi(\tau) = 1$. The terms $\varpi(\tau) = K(\tau)$ and $\boldsymbol{f}(\tau) = \mathbf{col}\,[1, g(\tau)]$ are in line with the definitions in [14, Theorem 1]. Finally, the inclusion of the IIs in [11] as special cases of (5) in Table I is demonstrated as follows.

Let $\alpha = 0$, $\beta = p \in \mathbb{N}_0$ with $\varpi(\tau) = (\tau - a)^p$ and $\boldsymbol{f}(\tau) = \boldsymbol{j}_d^{0,p}(\tau)_a^b$, for (5), where $\boldsymbol{j}_d^{0,p}(\tau)_a^b$ is defined (18b). Then we have

$$
\begin{aligned}
&\forall \boldsymbol{x}(\cdot) \in \mathscr{L}_{\varpi}^2(\mathcal{K};\mathbb{R}^n), \ \int_a^b (\tau-a)^p \boldsymbol{x}^{\top}(\tau) U \boldsymbol{x}(\tau)\mathsf{d}\tau \\
&\geq [*]\left(\mathsf{D}_{d,p} \otimes U\right) \int_a^b (\tau-a)^p \left(\boldsymbol{j}_d^{0,p}(\tau)_a^b \otimes I_n\right) \boldsymbol{x}(\tau)\mathsf{d}\tau \\
&= [*]\left(\mathsf{D}_{d,p} \otimes U\right) \Xi \int_a^b (\tau-a)^p \mathbf{Y}\left(\boldsymbol{\ell}_d(\tau) \otimes I_n\right) \boldsymbol{x}(\tau)\mathsf{d}\tau \\
&= [*]\left(\mathsf{D}_{d,p} \otimes U\right) \Xi \int_a^b \left(\boldsymbol{\ell}_{d+p}(\tau) \otimes I_n\right) \boldsymbol{x}(\tau)\mathsf{d}\tau \qquad (20)
\end{aligned}
$$

where $\mathsf{D}_{d,p}^{-1} = \bigoplus_{k=0}^{d} \frac{(b-a)^{p+1}}{2k+1+p}$ and

$$
\begin{aligned}
&\Xi := \mathbf{Y}(\mathsf{P} \otimes I_n), \ \ \mathbf{Y} := \left(\mathsf{J}_{0,p}\mathsf{L}^{-1}\right) \otimes I_n, \\
&\boldsymbol{j}_d^{0,p}(\tau)_a^b = \mathsf{J}_{0,p}\left[\tau^i\right]_{i=0}^d, \ \ \boldsymbol{\ell}_d(\tau) = \mathsf{L}\left[\tau^i\right]_{i=0}^d \\
&(\tau-a)^p \boldsymbol{\ell}_d(\tau) = \mathsf{P}\boldsymbol{\ell}_{d+p}(\tau).
\end{aligned} \qquad (21)
$$

with unique matrix parameters $\mathsf{J}_{0,p} \in \mathbb{R}_{[d+1]}^{(d+1)\times(d+1)}$ and $\mathsf{L} \in \mathbb{R}_{[d+1]}^{(d+1)\times(d+1)}$ and $\mathsf{P} \in \mathbb{R}^{(d+1)\times(d+p+1)}$. Using the multiplier $(b-a)^{-p}$, it shows that [11, eq.(27)] is equivalent to (20). Now using (20) with substitution $\dot{\boldsymbol{x}}(\tau) \to \boldsymbol{x}(\tau)$, we have

$$
\begin{aligned}
&\int_a^b (\tau-a)^p \dot{\boldsymbol{x}}^{\top}(\tau) U \dot{\boldsymbol{x}}(\tau)\mathsf{d}\tau \geq \\
&[*]\left[\Xi^{\top}(\mathsf{D}_d \otimes U)\,\Xi\right] \int_a^b \left(\boldsymbol{\ell}_{d+p}(\tau) \otimes I_n\right) \dot{\boldsymbol{x}}(\tau)\mathsf{d}\tau = \\
&\boldsymbol{\eta}_1^{\top}\Omega_1^{\top}\Xi^{\top}(\mathsf{D}_d \otimes U)\,\Xi\Omega_1\boldsymbol{\eta}_1 = \boldsymbol{\eta}_2^{\top}\Omega_2^{\top}\Xi^{\top}(\mathsf{D}_d \otimes U)\,\Xi\Omega_2\boldsymbol{\eta}_2 \qquad (22)
\end{aligned}
$$

where

$$
\begin{aligned}
&\boldsymbol{\eta}_1 := \begin{bmatrix} \boldsymbol{x}(b) \\ \boldsymbol{x}(a) \\ \int_a^b \left(\boldsymbol{\ell}_{d+p}(\tau) \otimes I_n\right) \boldsymbol{x}(\tau)\mathsf{d}\tau \end{bmatrix}, \ \ \boldsymbol{\eta}_2 := \begin{bmatrix} \boldsymbol{x}(b) \\ \boldsymbol{x}(a) \\ \int_a^b \left(\boldsymbol{\ell}_{d+p-1}(\tau) \otimes I_n\right) \boldsymbol{x}(\tau)\mathsf{d}\tau \end{bmatrix}, \\
&\Omega_1 = \left[\boldsymbol{\ell}_{d+p}(0) \ \ \boldsymbol{\ell}_{d+p}(-r) \ \ \Lambda_1\right] \otimes I_n, \\
&\Omega_2 = \left[\boldsymbol{\ell}_{d+p}(0) \ \ \boldsymbol{\ell}_{d+p}(-r) \ \ \Lambda_2\right] \otimes I_n
\end{aligned} \qquad (23)
$$

with $\Lambda_1 \in \mathbb{R}^{(d+p)\times(d+p)}$ and $\Lambda_2 \in \mathbb{R}^{(d+p)\times(d+p-1)}$ satisfying $\dot{\boldsymbol{\ell}}_{d+p}(\tau) = \Lambda_1\boldsymbol{\ell}_{d+p}(\tau) = \Lambda_2\boldsymbol{\ell}_{d+p-1}(\tau)$. Again by adjusting the factor $(b-a)^{-p}$ with (22), we show [11, eq.(34)] is a special case of (5).

Finally, the classical orthogonal polynomials [30] in the following table can also be used for $\boldsymbol{f}(\cdot)$ in (5) with appropriate weight functions, where $\lambda_d(\tau)_{\alpha}$ and $\hbar_d(\tau)$ represent Laguerre and Hermite polynomials, respectively, following their standard definitions.

| Names | Jacobi $j_d^{\alpha,\beta}(\tau)_{-1}^{1}$ | Laguerre $\lambda_d(\tau)_\alpha$ | Hermite $\hbar_d(\tau)$ |
|---|---|---|---|
| $\varpi(\tau)$ | $(1-\tau)^\alpha(1+\tau)^\beta$ | $x^\alpha \mathsf{e}^{-\tau}, \alpha > -1$ | $\mathsf{e}^{-\tau^2}$ |
| $\mathcal{K}$ | $[-1,1]$ | $[0,\infty)$ | $\mathbb{R}$ |
| $\|(\cdot)_d\|_\varpi^2$ | (14b) | $\frac{\Upsilon(d+\alpha+1)}{d!}$ | $2^d d! \sqrt{\pi}$ |

Classical Orthogonal Polynomials for $\boldsymbol{f}(\cdot)$

From the summary in Table I, it is clear that many IIs in the literature are mathematically equivalent. This could be attributed to the lack of awareness of the identities in (19), which connects the dots among many publications. The following proposition also includes a fundamentally important property that may have been overlooked in the literature.

**Proposition 1.** *Given the same functions and parameters in Theorem 1, we have*

$$\forall \boldsymbol{x}(\cdot) \in \mathscr{L}^2_\varpi(\mathcal{K};\mathbb{R}^n), \quad \int_{\mathcal{K}} \varpi(\tau)\boldsymbol{x}^\top(\tau)U\boldsymbol{x}(\tau)\mathsf{d}\tau$$
$$\geq \int_{\mathcal{K}} \varpi(\tau)\boldsymbol{x}^\top(\tau)F^\top(\tau)\mathsf{d}\tau\,(\mathsf{F}\otimes U)\int_{\mathcal{K}} \varpi(\tau)F(\tau)\boldsymbol{x}(\tau)\mathsf{d}\tau$$
$$= \int_{\mathcal{K}} \varpi(\tau)\boldsymbol{x}^\top(\tau)\Phi^\top(\tau)\mathsf{d}\tau\,(\mathsf{\Phi}\otimes U)\int_{\mathcal{K}} \varpi(\tau)\Phi(\tau)\boldsymbol{x}(\tau)\mathsf{d}\tau \quad (24)$$

*for all $G \in \mathbb{R}^{n\times n}_{[n]}$, where $\Phi(\tau) = \boldsymbol{\varphi}(\tau)\otimes I_n$ with $\boldsymbol{\varphi}(\tau) = G\boldsymbol{f}(\tau)$, and $\mathsf{\Phi}^{-1} := \int_{\mathcal{K}} \varpi(\tau)\boldsymbol{\varphi}(\tau)\boldsymbol{\varphi}^\top(\tau)\mathsf{d}\tau$ and $F(\tau)$, $\mathsf{F}$ are defined in the same way as in Theorem 1.*

*Proof.* Based on the definition of $\mathsf{\Phi}$, we see that

$$\mathsf{\Phi}^{-1} = \int_{\mathcal{K}} \varpi(\tau)\boldsymbol{\varphi}(\tau)\boldsymbol{\varphi}^\top(\tau)\mathsf{d}\tau = \int_{\mathcal{K}} \varpi(\tau)G\boldsymbol{f}(\tau)\boldsymbol{f}^\top(\tau)G^\top\mathsf{d}\tau$$
$$= G\int_{\mathcal{K}} \varpi(\tau)\boldsymbol{f}(\tau)\boldsymbol{f}^\top(\tau)\mathsf{d}\tau G^\top = G\mathsf{F}^{-1}G^\top \quad (25)$$

where $\mathsf{\Phi}^{-1}$ is well defined in respect that $G \in \mathbb{R}^{d\times d}_{[d]}$. With the help of (25) and the property in (1), we see that

$$\int_{\mathcal{K}} \varpi(\tau)\boldsymbol{x}^\top(\tau)\Phi^\top(\tau)\mathsf{d}\tau\;(\mathsf{\Phi}\otimes U)\int_{\mathcal{K}} \varpi(\tau)\Phi(\tau)\boldsymbol{x}(\tau)\mathsf{d}\tau$$
$$= \int_{\mathcal{K}} \boldsymbol{x}^\top(\tau)\left(\boldsymbol{\varphi}^\top(\tau)\otimes I_n\right)\mathsf{d}\tau\left[(G^{-1})^\top\mathsf{F}G^{-1}\otimes U\right]\times$$
$$\int_{\mathcal{K}} \left(\boldsymbol{\varphi}(\tau)\otimes I_n\right)\boldsymbol{x}(\tau)\mathsf{d}\tau = \int_{\mathcal{K}} \boldsymbol{x}^\top(t+\tau)\left(\boldsymbol{f}^\top(\tau)G^\top\otimes I_n\right)\times$$
$$\mathsf{d}\tau\left([*]\,(\mathsf{F}\otimes U)\left(G^{-1}\otimes I_n\right)\right)\int_{\mathcal{K}} \left(G\boldsymbol{f}(\tau)\otimes I_n\right)\boldsymbol{x}(t+\tau)\mathsf{d}\tau$$
$$= \int_{\mathcal{K}} \boldsymbol{x}^\top(t+\tau)F^\top(\tau)\mathsf{d}\tau\;(\mathsf{F}\otimes U)\int_{\mathcal{K}} F(\tau)\boldsymbol{x}(t+\tau)\mathsf{d}\tau,$$

which proves (24) by (5). ■

Proposition 1 reveals an important property that may have been consistently overlooked in the existing literature. Namely, using different linear combinations of $\{f_i(\cdot)\}_{i=1}^d$ does not change the LB value. This may explain why so many equivalent results in Table I have been developed, as many may have believed that using different $G$ for $\boldsymbol{f}(\cdot)$ could make the LBs tighter. A representative example is that using $\varpi(\tau) = 1$, $\boldsymbol{f}(\tau) = \boldsymbol{\ell}_d(\tau)$ with (5) produces exactly the same LB value as $\varpi(\tau) = 1$, $\boldsymbol{f}(\tau) = \left[\tau^i\right]_{i=0}^d$, based on Proposition 1 with $G = \mathsf{L}$ in (21). From a functional analysis perspective, Proposition 1 makes perfect sense because $\operatorname{span}_{i=0}^d \tau^i = \operatorname{span}_{i=0}^d \ell_i(\tau)$.

### *B. The Optimality of the LB of* (5)

Though the generality of (5) is evident because of the structures of $\boldsymbol{f}(\cdot) \in \mathscr{L}^2_\varpi(\mathcal{K};\mathbb{R}^d)$ and $\varpi(\cdot)$, it remains unclear what mathematical principle actually ensures the tightness of the LB in (5). The authors in [31] have shown that the optimality of the LB in (5) is guaranteed by the orthogonality principle, if $\boldsymbol{f}(\tau) = \boldsymbol{\ell}_d(\tau)$ with Legendre polynomials and $\varpi(\tau) = 1$. Here we can show that the optimality of the LB in (5) is always guaranteed by the principle of least squares approximations in Hilbert Space [21], even if the functions in $\boldsymbol{f}(\cdot)$ may not be orthogonal.

**Corollary 1.** *The II in* (5) *can be constructed via the least squares approximations in $\left(\mathscr{L}^2_\varpi(\mathcal{K};\mathbb{R})\,/\,\operatorname{Ker}(\|\cdot\|_\varpi), \langle\cdot,\cdot\rangle_\varpi\right)$ As a result, the optimality of* (5) *is guaranteed by the least squares principle. In other words, the LB in* (5) *is the best approximation of $\int_{\mathcal{K}} \varpi(\tau)\boldsymbol{x}^\top(\tau)U\boldsymbol{x}(\tau)\mathsf{d}\tau$ for all $\boldsymbol{x}(\cdot) \in \mathscr{L}^2_\varpi\left(\mathcal{K};\mathbb{R}^d\right)$ with given $U$, $\boldsymbol{f}(\cdot)$ and $\varpi(\cdot)$.*

*Proof.* Taking into account the definitions of matrix multiplication and the Kronecker product, $\boldsymbol{\varepsilon}(\tau) = \boldsymbol{x}(\tau) - F^\top(\tau)\boldsymbol{\omega}$ in (6) can be rewritten as

$$\boldsymbol{\varepsilon}(\tau) = \boldsymbol{x}(\tau) - (\boldsymbol{f}^\top(\tau)\otimes I_n)\boldsymbol{\omega} = \boldsymbol{x}(\tau) - \sum_{i=1}^d f_i(\tau)\boldsymbol{\omega}_i =$$
$$\boldsymbol{x}(\tau) - [\![\boldsymbol{\omega}_i]\!]_{i=1}^d\,\boldsymbol{f}(\tau) = \boldsymbol{x}(\tau) - \left[\boldsymbol{\kappa}_i^\top\right]_{i=1}^n\boldsymbol{f}(\tau) \quad (26)$$

where $\mathbb{R}^{dn} \ni \boldsymbol{\omega} := [\boldsymbol{\omega}_i]_{i=1}^d$ with $\boldsymbol{\omega}_i \in \mathbb{R}^n$, and $[\![\boldsymbol{\omega}_i]\!]_{i=1}^d = \left[\boldsymbol{\alpha}_j^\top\right]_{j=1}^n$ with $\boldsymbol{\alpha}_j \in \mathbb{R}^d$. Now by (26), it becomes clear that $\boldsymbol{\varepsilon}(\cdot)$ represents the vector containing a list of errors of approximation, and $\boldsymbol{\omega}$ contains all the approximation coefficient. Moreover, it shows that each row of $\boldsymbol{x}(\cdot) = [x_i(\tau)]_{i=1}^n$ is individually approximated by all the functions in $\boldsymbol{f}(\cdot) = [f_i(\cdot)]_{i=1}^d$.

Since $U \succ 0$, then $\exists N \in \mathbb{R}^{n\times n}_{[n]}$ such that $U = N^\top N$. With the semi-norm defined in (3), we seek to find $\boldsymbol{\lambda} = \arg\min_{\boldsymbol{\omega}\in\mathbb{R}^{dn}} \|N\boldsymbol{\varepsilon}(\cdot)\|_\varpi^2$ for a given $\boldsymbol{f}(\cdot) \in \mathscr{L}^2_\varpi(\mathcal{K};\mathbb{R}^d)$ that minimizes $\|N\boldsymbol{\varepsilon}(\cdot)\|_\varpi^2 = \int_{\mathcal{K}} \varpi(\tau)\boldsymbol{\varepsilon}^\top(\tau)U\boldsymbol{\varepsilon}(\tau)\mathsf{d}\tau$. As we can always write $U = N^\top N \succ 0$, then there exist $c_1; c_2 > 0$, $c_1\|\boldsymbol{\varepsilon}(\cdot)\|_\varpi^2 \leq \|N\boldsymbol{\varepsilon}(\cdot)\|_\varpi^2 \leq c_2\|\boldsymbol{\varepsilon}(\cdot)\|_\varpi^2$. Hence, it follows that $\arg\min_{\boldsymbol{\omega}\in\mathbb{R}^{dn}} \|N\boldsymbol{\varepsilon}(\cdot)\|_\varpi^2 = \arg\min_{\boldsymbol{\omega}\in\mathbb{R}^{dn}} \|\boldsymbol{\varepsilon}(\cdot)\|_\varpi^2$. Now coefficient $\boldsymbol{\lambda}$ can be obtained by making use of the least squares principle [21] in the Hilbert Space

$$\left(\mathscr{L}^2_\varpi(\mathcal{K};\mathbb{R})\,/\,\operatorname{Ker}(\|\cdot\|_\varpi), \langle\cdot,\cdot\rangle_\varpi\right) \quad (27)$$

with inner product $\langle f,g\rangle_\varpi = \int_{\mathcal{K}} \varpi(\tau)f(\tau)g(\tau)\mathsf{d}\tau$. Namely, for any $x_i(\cdot) \in \mathscr{L}^2_\varpi(\mathcal{K};\mathbb{R})$ we want to find the closest point $\widehat{f_i}(\cdot) = \boldsymbol{\kappa}_i^\top\boldsymbol{f}(\tau)$ to $x_i(\cdot)$ with the optimal coefficient $\boldsymbol{\lambda}$, where $\widehat{f_i}(\cdot) \in \operatorname{span}_{i=1}^d f_i(\cdot)$ and $\boldsymbol{\lambda} = [\![\boldsymbol{\lambda}_i]\!]_{i=1}^d = \left[\boldsymbol{\kappa}_j^\top\right]_{j=1}^n$. For more expositions on this topic, see [21, section 10.2].

Based on the least squares principle on [21, page 182], $\boldsymbol{\lambda} = \arg\min_{\boldsymbol{\omega}\in\mathbb{R}^{dn}} \|\boldsymbol{\varepsilon}(\cdot)\|_\varpi^2$ is unique and satisfy the relation

$$\mathsf{F}^{-1}\boldsymbol{\kappa}_i = \int_{\mathcal{K}} \varpi(\tau)x_i(\tau)\boldsymbol{f}(\tau)\mathsf{d}\tau \in \mathbb{R}^d, \quad i = 1,\ldots,n$$

which can then be denoted by the compact form

$$\mathsf{F}^{-1}\left[\!\left[\boldsymbol{\kappa}_i\right]\!\right]_{i=1}^{n}=\mathsf{F}^{-1}\left[\boldsymbol{\lambda}_i^{\top}\right]_{i=1}^{d}=\int_{\mathcal{K}}\varpi(\tau)\boldsymbol{f}(\tau)\boldsymbol{x}^{\top}(\tau)\mathsf{d}\tau\in\mathbb{R}^{d\times n}.$$

Thus we conclude $[\![\boldsymbol{\lambda}_i]\!]_{i=1}^{d}=\int_{\mathcal{K}}\varpi(\tau)\boldsymbol{x}(\tau)\boldsymbol{f}^{\top}(\tau)\mathsf{d}\tau\mathsf{F}$. Now by the definition of matrix multiplications and a new function $\widetilde{\boldsymbol{f}}(\tau):=\left[\widetilde{f}(\tau)\right]_{i=1}^{d}:=\mathsf{F}\boldsymbol{f}(\tau)$, we see that $[\![\boldsymbol{\lambda}_i]\!]_{i=1}^{d}=\left[\!\left[\int_{\mathcal{K}}\varpi(\tau)\boldsymbol{x}(\tau)\widetilde{f}_i(\tau)\mathsf{d}\tau\right]\!\right]_{i=1}^{d}$, which indicates

$$\boldsymbol{\lambda}_i=\int_{\mathcal{K}}\varpi(\tau)\boldsymbol{x}(\tau)\widetilde{f}_i(\tau)\mathsf{d}\tau,\ \ \forall i=1,\cdots,d. \tag{28a}$$

Utilizing $\boldsymbol{\lambda}=[\boldsymbol{\lambda}_i]_{i=1}^{d}\in\mathbb{R}^{dn}$ with (28a), we see that

$$\begin{aligned}\boldsymbol{\lambda}&=\left[\int_{\mathcal{K}}\varpi(\tau)\widetilde{f}_i(\tau)\boldsymbol{x}(\tau)\mathsf{d}\tau\right]_{i=1}^{d}=\int_{\mathcal{K}}\varpi(\tau)(\widetilde{\boldsymbol{f}}(\tau)\otimes I_n)\boldsymbol{x}(\tau)\mathsf{d}\tau\\&=(\mathsf{F}\otimes I_n)\int_{\mathcal{K}}\varpi(\tau)(\boldsymbol{f}(\tau)\otimes I_n)\boldsymbol{x}(\tau)\mathsf{d}\tau=(\mathsf{F}\otimes I_n)\boldsymbol{\vartheta}\end{aligned}\tag{28b}$$

with $\widetilde{\boldsymbol{f}}(\tau)=\mathsf{F}\boldsymbol{f}(\tau)$. Since (12) is also obtained with $\boldsymbol{\omega}=\boldsymbol{\lambda}=(\mathsf{F}\otimes I_n)\boldsymbol{\vartheta}$ corresponding to the largest LB for (11), this proves that the optimality of the LB in (5) is guaranteed by the least squares principle [21]. ■

*Remark* 3. Corollary 1 establishes that the optimality derived from the orthogonality (least squares) principle is preserved for any $\boldsymbol{f}(\cdot)\in\mathscr{L}^2_{\varpi}(\mathcal{K};\mathbb{R}^d)$ satisfying (4), even if the functions in $\boldsymbol{f}(\cdot)$ are not mutually orthogonal. This is a significant conclusion as it challenges the common belief that the optimality of the orthogonality (least squares) principle [21] can only be obtained with orthogonal functions.

Since the LB in (5) is the best approximation of the upper bound integral with given $U$, $\boldsymbol{f}(\cdot)$ and $\varpi(\cdot)$, we show in the following theorem that the LB value increases if the dimension of $\boldsymbol{f}_d(\tau)$ is enlarged.

**Theorem 2.** *Let $\varpi(\cdot)$, $\mathcal{K}$ in Theorem 1 be given, and assume $\boldsymbol{f}_{d+1}(\cdot)=[f_i(\cdot)]_{i=1}^{d+1}\in\mathscr{L}^2_{\varpi}(\mathcal{K};\mathbb{R}^{d+1})$ satisfying*

$$\mathsf{F}_{d+1}^{-1}=\int_{-r}^{0}\varpi(\tau)\boldsymbol{f}_{d+1}(\tau)\boldsymbol{f}_{d+1}^{\top}(\tau)\mathsf{d}\tau\succ 0 \tag{29}$$

*with reference to* (4)*. Then for all $d\in\mathbb{N}$,*

$$\begin{aligned}&\int_{\mathcal{K}}\varpi(\tau)\boldsymbol{x}^{\top}(\tau)U\boldsymbol{x}(\tau)\mathsf{d}\tau\\&\quad\geq[*]\left(\mathsf{F}_{d+1}\otimes U\right)\int_{\mathcal{K}}\varpi(\tau)F_{d+1}(\tau)\boldsymbol{x}(\tau)\mathsf{d}\tau\\&\qquad\geq[*]\left(\mathsf{F}_{d}\otimes U\right)\int_{\mathcal{K}}\varpi(\tau)F_{d}(\tau)\boldsymbol{x}(\tau)\mathsf{d}\tau\end{aligned}\tag{30}$$

*where $\mathsf{F}_d^{-1}=\int_{\mathcal{K}}\varpi(\tau)\boldsymbol{f}_d(\tau)\boldsymbol{f}_d^{\top}(\tau)\mathsf{d}\tau\succ 0$ and $F_d(\tau)=\boldsymbol{f}_d(\tau)\otimes I_n$ for all $d\in\mathbb{N}$.*

*Proof.* The proof is omitted due to limited space. It will be presented in the journal version. ■

We set forth the above result as a theorem because it reveals one of the most attractive features of the proposed II in (5). Namely, we can increase the LB value of (5) simply by adding more linearly independent functions to $\boldsymbol{f}(\cdot)$. A special case of this property was first remarked in [1] where the proof was achieved based on $\boldsymbol{f}_d(\tau)=\boldsymbol{\ell}_d(\tau)$ with the help of the orthogonality of Legendre polynomials. Our Theorem 2, however, only requires that $\boldsymbol{f}_d(\cdot)$ contains linearly independent functions. Hence, not only Theorem 2 does presents a significant generalization of the idea in [1], but also reveals that the successive improvement of the LB value does not require $\{f_i(\cdot)\}_{i=1}^{d}$ being mutually orthogonal.

*Remark* 4. Another way to understand (30) is to realize

$$\forall d\in\mathbb{N},\ \operatorname{span}_{i=1}^{d}f_i(\cdot)\subset\operatorname{span}_{i=1}^{d+1}f_i(\cdot)\subset\mathscr{L}^2_{\varpi}\left(\mathcal{K};\mathbb{R}^d\right), \tag{31}$$

since (29) indicates that $\{f_i(\cdot)\}_{i=1}^{d+1}$ are linearly independent. The inclusions in (31) provide an intuitive explanation for why adding more linearly independent functions to $\boldsymbol{f}(\cdot)$ can increase the LB value of (5).

Given the result in Theorem 2, a natural question arises: to what extent can the LB value be increased? The following corollary provides an interesting answer to this question, showing that the LB in (5) can asymptotically approach the integral upper bound with $\mathcal{K}=[a,b]$ and appropriate functions $\{f_i(\cdot)\}_{i=1}^{\infty}$ as $d\to\infty$.

**Corollary 2.** *Assume that the Hilbert space in* (27) *is separable. Let $\varpi(\cdot)$ and $U\succ 0$ be given, and function series $\{f_i(\cdot)\}_{i=1}^{\infty}$ is a Schauder basis of* (27)*. Then*

$$\begin{aligned}&\int_{\mathcal{K}}\varpi(\tau)\boldsymbol{x}^{\top}(\tau)U\boldsymbol{x}(\tau)\mathsf{d}\tau=\\&\qquad\lim_{d\to\infty}[*]\left(\mathsf{F}_d\otimes U\right)\int_{\mathcal{K}}\varpi(\tau)F_d(\tau)\boldsymbol{x}(\tau)\mathsf{d}\tau.\end{aligned}\tag{32}$$

*Proof.* The proof is omitted due to limited space. It will be presented in the journal version. ■

If $\{f_i(\cdot)\}_{i=1}^{\infty}$ is an orthonormal basis of the Hilbert space in (27), then (32) is the common orthonormal decomposition utilized in [1], [20] as a special case of the Bessel inequality [21, Theorem 10.31]. This shows Corollary 2 generalizes the framework developed in [1], [20].

With Corollary 2, the result in Theorem 2 becomes even more powerful, as it shows that the LB in (5) can converge exactly to the upper bound integral with a Schauder basis $\{f_i(\cdot)\}_{i=1}^{\infty}$ as $d\to\infty$. Since a Schauder basis may not be orthogonal nor an orthonormal basis, it precludes the need of using the Bessel inequality [21, Theorem 10.31] for the analysis of (5) as in [1], [20]. In consequence, we have demonstrated via Corollary 2 that the valuable property of asymptotic convergence does not require $\{f_i(\cdot)\}_{i=1}^{\infty}$ to be an orthonormal basis.

## III. Another Integral Inequality

The construction of (5) is more of a direct approach as $U$ is the only matrix that can be interpreted as a decision variable in an optimization context. With the belief that having slack matrix variables in IIs may reduce the conservatism of their LBs, researchers have proposed numerous IIs [15]–[17] whose LBs have more matrix variables than the original $U$. These IIs are referred to as the free matrix type (FMT), and have been utilized in situations where $\mathcal{K}$ contains parameters that are functions of other variables, such as the case of $\mathcal{K} = [-r(t), 0]$ with a time-varying delay function $r(t)$.

We propose a general FMTII in this section that allows us to utilize the general $\boldsymbol{f}(\cdot)$ function as in (5). For this reason, our FMTII can generalize all existing IIs that have a free-matrix structure. Moreover, an interesting and crucial result shows that the maximum value of the LB of the proposed FMTII is equivalent to the one in (5), under the same $\boldsymbol{f}(\cdot)$, $\varpi(\cdot)$, and $U$.

Due to limited space, the content of this section is omitted in this note, which will be presented in the journal version. As a reference, the proposed FMTII has the structure

$$\begin{aligned}\int_{\mathcal{K}} \varpi(\tau)\boldsymbol{x}^\top(\tau)U\boldsymbol{x}(\tau)\mathrm{d}\tau &\geq \mathsf{Sy}\left(\boldsymbol{\vartheta}^\top \widehat{X}\mathbf{z}\right) - \mathbf{z}^\top \mathsf{W}\mathbf{z} \\ \int_{\mathcal{K}} \varpi(\tau)\boldsymbol{x}^\top(\tau)U\boldsymbol{x}(\tau)\mathrm{d}\tau &\geq \mathbf{y}^\top \left[\mathsf{Sy}\left(\Upsilon^\top \widehat{X}\right) - \mathsf{W}\right]\mathbf{y}\end{aligned} \tag{33}$$

for all $\mathbf{z} \in \mathbb{R}^{\rho n}$ and $\boldsymbol{x}(\cdot) \in \mathscr{L}^2_{\varpi}(\mathcal{K};\mathbb{R}^n)$, where $\boldsymbol{\vartheta} = \int_{\mathcal{K}} \varpi(\tau)F(\tau)\boldsymbol{x}(\tau)\mathrm{d}\tau$, $\widehat{X} = \mathbf{Col}_{i=1}^d X_i \in \mathbb{R}^{dn\times\rho n}$, and $\mathsf{W} = \int_{\mathcal{K}} \varpi(\tau)(\boldsymbol{f}^\top(\tau)\otimes I_{\rho n})Y(\boldsymbol{f}(\tau)\otimes I_{\rho n})\mathrm{d}\tau \in \mathbb{S}^{\rho n}$, and $\Upsilon \in \mathbb{R}^{dn\times\rho n}$, $\mathbf{y} \in \mathbb{R}^{\rho n}$ satisfy $\forall \boldsymbol{x}(\cdot) \in \mathscr{L}^2_{\varpi}(\mathcal{K};\mathbb{R}^n)$, $\Upsilon\mathbf{y} = \boldsymbol{\vartheta}$, with some $Y \in \mathbb{S}^{\rho d n}$ and $X = [\![X_i]\!]_{i=1}^d \in \mathbb{R}^{n\times\rho d n}$ satisfying $\begin{bmatrix} U & -X \\ * & Y\end{bmatrix} \succ 0$

## IV. Conclusion

In this note, some general classes of IIs have been introduced, where the integral kernels of the LBs can be any function in $\mathscr{L}^2_{\varpi}(\mathcal{K};\mathbb{R})$ satisfying (4) with a general weight function $\varpi(\tau)$ in (3). It has been shown that the proposed IIs generalize many existing examples in the literature. We have demonstrated that the optimality of (5) with $U \succ 0$ is guaranteed by the least squares principle [21] as summarized in Corollary 1. Moreover, the equivalence between the LBs of many existing IIs in Table I are established using Proposition 1 together with the generality of (5). Additional contents related to our FMTII with numerous interesting conclusions will be presented in the journal version of this paper. The proposed IIs have great potential to be applied in wider contexts such as the stability analysis of PDE-related systems, time-delay systems, or other types of distributed parameters whenever the situations are suitable.